\documentclass[12pt]{amsart}
\usepackage{psfrag}
\usepackage[usenames]{color}
\usepackage[normalem]{ulem}%\uline,\uuline,\uwave,\sout,\xout
\usepackage{citesort}

\DeclareMathOperator{\lcm}{lcm} 
\newcommand{\orb}{\mathcal O}
\newcommand{\clorb}{\overline{\mathcal O}}
\newcommand{\R}{\mathbb R}
 
\newcommand{\Z}{\mathbb Z}
\newcommand{\T}{\mathbb T}

\newcommand{\Q}{\mathbb Q}
\newcommand{\bs}{\backslash} 
 
\DeclareMathOperator{\cl}{cl}

\newcommand{\ep}{\varepsilon}

\theoremstyle{plain} 
\newtheorem{thm}{Theorem}
\newtheorem{cor}[thm]{Corollary} 
\newtheorem{prop}[thm]{Proposition}
\newtheorem{lemma}[thm]{Lemma}

\theoremstyle{definition}

\theoremstyle{remark} 
\newtheorem{remark}[thm]{Remark}

\begin{document} 
    
\title[Itineraries of rigid rotations]{Itineraries of rigid rotations and diffeomorphisms of the circle}

\author{David Richeson} \author{Paul Winkler}  \author{Jim Wiseman} \address{Dickinson
College\\ Carlisle, PA 17013} \email{richesod@dickinson.edu} \address{Dickinson
College\\ Carlisle, PA 17013}\email{winklerp@dickinson.edu}
\address{Agnes Scott College \\ Decatur, GA 30030}
\email{jwiseman@agnesscott.edu}

\date{\today}
\keywords{symbolic dynamics, rotations}
\subjclass[2000]{}
\begin{abstract}
We examine the itinerary of $0\in S^{1}=\R/\Z$ under the rotation by $\alpha\in\R\bs\Q$.  The motivating question is: if we are given only the itinerary of $0$ relative to $I\subset S^{1}$, a finite union of closed intervals, can we recover $\alpha$ and $I$?  We prove that the itineraries do determine $\alpha$ and $I$ up to certain equivalences.  Then we present elementary methods for finding $\alpha$ and $I$.  Moreover, if $g:S^{1}\to S^{1}$ is a $C^{2}$, orientation preserving diffeomorphism with an irrational rotation number, then we can use the orbit itinerary to recover the rotation number up to certain equivalences.
\end{abstract}

\maketitle
A useful and common technique for analyzing discrete dynamical systems is to partition the space and study the itineraries and corresponding shift spaces of orbits.  Often the dynamics of the original map is complicated (or chaotic) and the shift space provides a convenient way of extracting properties of the original dynamical system.   In this paper, however, we consider orbit itineraries for one of the most elementary and well-understood dynamical systems---the rigid rotation of a circle.  

Let $r_{\alpha}:S^{1}=\R/\Z\to S^{1}$ be the rotation $r_{\alpha}(z)=z+\alpha \mod 1$, where $\alpha\in\R$.  Let $I$ be a finite union of closed intervals (with nonempty interiors) in $S^{1}$ that is neither $\emptyset$ nor $S^{1}$.  The {\em itinerary} for $0\in S^{1}$ is $(a_{0},a_{1},a_{2},\ldots)$, where \[a_{i}=\begin{cases}0 & \text{if } r_{\alpha}^{i}(0)\not\in I\\1 & \text{if } r_{\alpha}^{i}(0)\in I \end{cases}\]

For example, when $\alpha=\frac{\sqrt{3}}{15}$ and $I=[0,\frac{1}{4}]$ the itinerary of $0$ is \[(1,1,1,0,0,0,0,0,0,1,1,0,0,0,0,0,0,0,1,1,0,0,0,0,0,0,1,1,1,\]\[0,0,0,0,0,0,1,1,0,0,0,0,0,0,0,1,1,0,0,0,0,0,0,1,1,1,\ldots)\]
\label{sampleitinerary}

In this paper we investigate these orbit itineraries.  In particular we answer the following question: given an itinerary and no other information, can we find $\alpha$ and $I$?  We will see that if $\alpha$ is  irrational  and $I$ has no rotational symmetries in $S^{1}$, then we can find the fractional part of $\alpha$ up to sign, and for these two $\alpha$-values we can find the corresponding $I$.  If $I$ has rotational symmetries then (assuming we do not know the order of this symmetry) we can find $\alpha$ or $-\alpha$ up to an integer multiple.  We can say much less when $\alpha$ is rational, so in this paper we will focus on irrational rotations.  (In the case that $I$ is a single interval and $\alpha$ is rational, $\alpha$ and $I$ can in general be computed only up to a certain level of accuracy; in many cases, $\alpha$ can be computed exactly.  See \cite{AGM,STZ} for more details and discussion.)

The paper is organized as follows.  We define notation and discuss previous work in Section~\ref{sec:notaPrev}.  In Section~\ref{sec:unique}, we show that, except for certain symmetries, itineraries are unique (that is, if two itineraries are the same, then both the angles and the intervals must also be the same, up to symmetry, for both).  In Section~\ref{sec:findI}, we give an easy method for finding $I$ given $\alpha$. In Section~\ref{sec:singleInt} we discuss the more difficult problem of finding $\alpha$.  We give a method that works well in the case that $I$ is a single interval but is not certain to give a good estimate (using only a finite portion of the itinerary) in the case that $I$ comprises multiple intervals.  Finally, in Section~\ref{sec:diffeos} we apply the results for rotations to orientation preserving diffeomorphisms of the circle.

We are grateful to the referees for their valuable suggestions for improving the paper, and to Alan Koch for helpful conversations.

\section{Notation and previous work}
\label{sec:notaPrev}

Again, let $r_{\alpha}:S^{1}=\R/\Z\to S^{1}$ be the rotation $r_{\alpha}(z)=z+\alpha \mod 1$, where $\alpha\in\R$, and let $I$ be a finite union of closed intervals (with nonempty interiors) in $S^{1}$ that is neither $\emptyset$ nor $S^{1}$.
We should point out that the choice of whether to include or exclude the endpoints of our intervals in $I$ was arbitrary.  Unless an endpoint of $I$ is a multiple of $\alpha$, the orbit of $0$ will not include it.  At most, the infinite orbit will land on the endpoints only finitely many times.  For convenience, we assume that $I$ is closed.  We let $l(I)$ denote the sum of the lengths of the intervals in $I$.

To simplify notation let $\{y\}=y-\lfloor y\rfloor$ denote the fractional part of $y$, and let $[y]=\min\{\{y\},\{-y\}\}$.  For example, $\{-1.3\}=0.7$ and $[-1.3]=0.3$.  We will also repeatedly be lazy with notation and write $z$ when we mean $z\in\R$ (mod 1), $z\in[0,1)$, and $z\in S^{1}$. Let $\orb(z)=\bigcup_{i=0}^\infty r_{\alpha}^i(z)$ and $\clorb(z)=\cl(\orb(z))$ be the orbit and orbit closure of $z$, respectively.

The dynamics of rigid rotations has been studied extensively.  The literature is broad and diverse, and we only touch on it here.  Hedlund studied the itinerary of the point $0\in S^{1}$ under a rigid rotation $\alpha$ and interval $I=[0,\alpha)$ (\cite{Hed}).  The associated shift space, called the {\em Sturmian shift}, is well-studied and has many interesting properties.  According to Coven and Nitecki (\cite{CN}) this was the first example of symbolic dynamics.  Slater and others looked at the more general case of the itinerary of $0\in S^{1}$ under a rotation by $\alpha$ and interval $I=[0,\beta)$ (\cite{S1,S2,AB}).  He proved the so-called ``three gap theorem,'' that implies that maximal blocks of the form ``$0,0,\ldots,0$'' can have at most three possible lengths, and if the lengths are $a$, $b$, and $c$ (with $c$ the largest), then $c=a+b+1$ (and the same is true for the maximal $1,1,\ldots,1$ blocks).  For example, in the itinerary given earlier, the $0,0,\ldots,0$ blocks have lengths 0 (since there are two adjacent $1$'s), 6, and 7 (likewise, the $1,\ldots,1$ blocks have length 0, 2, and 3). These itineraries have also been studied by decomposing them into Sturmian itineraries (\cite{BV,Did1,Hub}), with continued-fraction-like approaches (\cite{Ada,Did2}), and using interval exchange transformations (\cite{BCZ,FHZ,LN}).  Siegel, Tresser and Zettler investigated the still more general case of itineraries for homeomorphisms with rotation number $\alpha$ and corresponding interval $I=[0,\beta)$ (\cite{STZ}).  For a more complete bibliography and background, see \cite{Ada} and \cite{STZ}.

In this paper we take a more elementary approach than in most of the previous investigations.

\section{Uniqueness of itineraries}
\label{sec:unique}

We would like to show that for irrational rotations, the itineraries uniquely determine the rotation and the set $I$.  More specifically, suppose  $\alpha_{1},\alpha_{2}\in[0,1)\bs\Q$ and $I_{1},I_{2}\subset S^{1}$ are finite unions of closed intervals.  We would like to show that $(\alpha_{1},I_{1})=(\alpha_{2},I_{2})$ if and only if the associated itineraries are the same.  Clearly that is not true.  We encounter non-uniqueness corresponding to the symmetry of clockwise or counterclockwise rotations.  The itinerary for $\alpha_{1}=\sqrt{2}-1$, $I_{1}=[1/4,1/2]$ is the same as the itinerary for $\alpha_{2}=2-\sqrt{2}$, $I_{2}=[1/2,3/4]$.  Rotational symmetries for $I$ can also cause non-uniqueness.  We say that an interval $I\subset \R/\Z$ has \emph{$n$-fold rotational symmetry} if $I+1/n=I\subset\R/\Z$. These itineraries are also the same as the itineraries for $\alpha_{3}=(\sqrt{2}-1)/2$, $I_{3}=[1/8,1/4]\cup[5/8,3/4]$, and $\alpha_{4}=(\sqrt{2}-1)/3$, $I_{4}=[1/12,1/6]\cup [5/12,1/2]\cup [3/4,5/6]$.  However, these are the only causes of non-uniqueness.  We will prove the following theorem, which is a somewhat stronger version of \cite[Lemma~4.1]{Did1}. 

\begin{thm}\label{thm:unique}
Let $\alpha_{1},\alpha_{2}\in[0,1/2]\bs \Q$ and let $I_{1},I_{2}\subset S^{1}$ be finite unions of closed intervals with no rotational symmetries.  Then the corresponding itineraries are the same if and only if $\alpha_{1}=\alpha_{2}$ and $I_{1}=I_{2}$.
\end{thm}

If we allow $\alpha_{1}$ and $\alpha_{2}$ to be any irrational numbers we may rephrase the theorem as follows.

\begin{cor}
Let $\alpha_{1},\alpha_{2}\in\R\bs \Q$ and let $I_{1},I_{2}\subset S^{1}$ be finite unions of closed intervals with no rotational symmetries.  Then the corresponding itineraries are the same if and only if $\{\alpha_{1}\}=\{\alpha_{2}\}$ and $I_{1}=I_{2}$ or $\{\alpha_{1}\}=1-\{\alpha_{2}\}$ and $I_{1}=-I_{2}$.
\end{cor}

%\begin{proof}
%Now, suppose $[\alpha_{1}]\ne[\alpha_{2}]$.  It suffices to show that there is some $k>0$ such that $f_{\alpha_{1}}^{k}(z)\in I$ and $f_{\alpha_{2}}^{k}(z)\not\in I$.  Equivalently, if we let $F:S^{1}\times S^{1}\to S^{1}\times S^{1}$ be the map $F(x,y)=(f_{\alpha_{1}}(x),f_{\alpha_{2}}(y))=(x+\alpha_{1},y+\alpha_{2}) \mod 1$, then we must show that for some $k>0$,  $F^{k}(z,z)\in D$, where $D=([0,\beta)\times[\beta,1))\cup([\beta,1)\times[0,\beta))$.   
%    
%First consider the case that $\alpha_{2}/\alpha_{1}$ is rational.  The orbit of $(z,z)$ lies on the line $y-z=(\alpha_{2}/\alpha_{1})(x-z)$ in the torus. Since the slope of the line is rational, it is an invariant circle.  We know that $(z,z)$ is not periodic, so its orbit must be dense in the circle.  The only circle that avoids $I\times (S^{1}\bs I)$ is the diagonal, which this circle is not.  Thus $F^{k}(z,z)\in I\times (S^{1}\bs I)$ for some $k>0$
%    
%Now suppose that $\alpha_{2}/\alpha_{1}$ is irrational.  Then $F$ is uniquely ergodic (see \cite[Sec. 1.4]{KH}).  In particular, the orbit of $(z,z)$ is dense in the torus.  Thus, $F^{k}(z,z)\in I\times (S^{1}\bs I)$ for some $k>0$.
%\end{proof}

\begin{figure}[ht]
\centering
\psfrag{0}{$0$}
\psfrag{1}{$1$}
\psfrag{i1}{$I_{1}$}
\psfrag{i2}{$I_{2}$}
\psfrag{b}{$B$}
\includegraphics{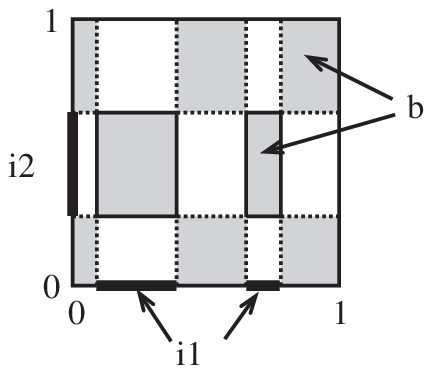}
\caption{The set $B=(I_{1}\times I_{2})\cup(I_{1}^{c}\times I_{2}^{c})$}
\label{fig:beta1}
\end{figure}

We can rephrase the uniqueness question as one about translations on the torus.  Consider the map $R=R_{\alpha_{1},\alpha_{2}}=r_{\alpha_{1}}\times r_{\alpha_{2}}:S^{1}\times S^{1}\to S^{1}\times S^{1}$ given by $R(z_{1},z_{2})=(z_{1}+\alpha_{1},z_{2}+\alpha_{2})$ and let $B=(I_{1}\times I_{2})\cup(I_{1}^{c}\times I_{2}^{c}),$ as in Figure \ref{fig:beta1} (where $I_{j}^{c}=S^{1}\bs I_{j}$). Then the itineraries corresponding to $(\alpha_{1},I_{1})$ and $(\alpha_{2},I_{2})$ are the same if and only if the orbit of $(0,0)$  lies in $B$.  Thus, in order to prove Theorem \ref{thm:unique} we must investigate the properties of the orbit of $(0,0)$ under $R_{\alpha_{1},\alpha_{2}}$.  

We say that $\alpha_{1}$, $\alpha_{2}$, and 1 are \emph{rationally related} if there exist integers $a$, $b$, and $c$, not all zero, such that $a\alpha_{1}+b=c\alpha_{2}$ (that is, $\alpha_{1}$, $\alpha_{2}$, and 1 are linearly dependent in the vector space $\R$ over $\Q$). The following theorem describes the dynamics of $R$.

\begin{prop}
\label{prop:torusorbits}
Suppose $\alpha_{1},\alpha_{2}\in\R$ and $R:\T^{2}\to\T^{2}$ is given by $R(x,y)=(r_{\alpha_{1}}(x),r_{\alpha_{2}}(y))$. 
\begin{enumerate}
\item
If $\alpha_{1}$ and $\alpha_{2}$ are rational, then every point is periodic.
\item
If $\alpha_{1}$ and $\alpha_{2}$ are rationally related but not both rational, then the closure of any orbit is a finite collection of circles.  
\item
If $\alpha_{1}$ and $\alpha_{2}$ are not rationally related, then the closure of every orbit is $\T^{2}$.
\end{enumerate} 
Moreover, $R$ restricted to the closure of any orbit is uniquely ergodic (that is, it has only one invariant Borel probability measure).
\end{prop}
\begin{proof}
(1) If $\alpha_{1}=p_{1}/q_{1}$ and $\alpha_{2}=p_{2}/q_{2}$ are reduced fractions, then every point is periodic with period $\lcm(|q_{1}|,|q_{2}|)$.

(2) This is essentially Exercise 1.4.1 in \cite{KH}.  We will prove this for the orbit of $(0,0)$; the proof for other points is similar.  Suppose $a\alpha_{1}+b=c\alpha_{2}$ for some $a,b,c\in \Z$ where $a$, $b$, and $c$ are not all zero and have no common factors.  Either $a$ or $c$ is nonzero, so without loss of generality, assume that $c\ne 0$.  It is clear that the point $(m\alpha_{1},m\alpha_{2})$ lies on the line $y=(a/c)x+mb/c$ in $\R^{2}$ and that the set of all such lines projects to a finite number of circles in $\T^{2}$.  Thus $\orb(0,0)$ is contained in a finite number of circles ($g$, say). $R$ permutes these circles cyclically and $R^{g}$ restricted to one of these circles is conjugate to an irrational rotation; thus $\clorb(0,0)$ is the set of $g$ circles.

(3)  In this case, $R$ is minimal, that is, every orbit is dense (see \cite[\S1.4]{KH}).

Unique ergodicity is clear if the orbit is periodic, and well known for irrational rotations and translations (see \cite[\S4.2]{KH}).

\end{proof}
\begin{remark}
In fact, we can say exactly how many circles there are in case (2) of Proposition \ref{prop:torusorbits}. If one of the $\alpha_{i}$ is rational, say $\alpha_{i}=p/q$ (reduced), then $\clorb(0,0)$ is $|q|$ circles.  Suppose neither $\alpha_{i}$ is rational and $a\alpha_{1}+b=c\alpha_{2}$ for some $a,b,c\in \Z$ where $a$, $b$, and $c$ are not all zero and have no common factors.  Let $g=\gcd(|a|,|c|)$. Observe that the lines $y=(a/c)x+mb/c$ ($m\in\Z$) project to $n$ circles in $\T^{2}$ if and only if $n$ is the smallest positive integer for which $y=(a/c)x+nb/c$ is the same circle as $y=(a/c)x$ in $\T^{2}$.  In other words, $n$ is the smallest positive integer for which $y=(a/c)x+nb/c$ passes through a point $(x,y)\in\Z^{2}$.  It is a basic fact from algebra (see Theorem 0.2 in \cite{Gal}, for instance) that $cy-ax=nb$ has integer solutions if and only if $nb$ is a multiple of $g$.  If $b=0$, then $n=1$ is the smallest positive value of $n$ for which this is true, so $\clorb(0,0)$ is a single circle.  If $b\ne 0$, then because $\gcd(b,g)=1$, the smallest positive value of $n$ for which this is true is $n=g$, so $\clorb(0,0)$ is $g$ circles.
\end{remark}

\begin{figure}[ht]
\centering
\psfrag{y1}{$y=x/2$}
\psfrag{y2}{$y=x/2+1/4$}
\includegraphics{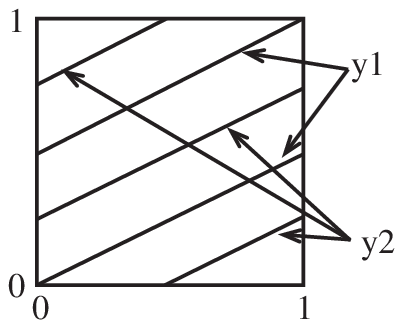}
\caption{The set $\clorb(0,0)$ for $\alpha_{1}=(\sqrt{2}-2)/2$ and $\alpha_{2}=(\sqrt{2}+3)/4$}
\label{fig:circlesontorus}
\end{figure}

For example, $\alpha_{1}=(\sqrt{2}-2)/2$ and $\alpha_{2}=(\sqrt{2}+3)/4$ are rationally related, with $2\alpha_{1}+5=4\alpha_{2}$.  The closure of the orbit of $(0,0)$ is two circles, namely $y=x/2$ and $y=x/2+5/4$, or equivalently on the torus, $y=x/2+1/4$  (see Figure \ref{fig:circlesontorus}).

\begin{prop}\label{prop:ratdep}
Suppose $I_{1},I_{2}\subset S^{1}$ are finite unions of closed intervals and $\alpha_{1},\alpha_{2}\in\R\bs \Q$ are rationally related (with $a$, $b$, $c$, and $g$ as above).  The corresponding itineraries are the same if and only if $I_{1}$ has $|a|$-fold symmetry, $I_{2}$ has $|c|$-fold symmetry, and $I_{2}=\{ax/c+mb/c:x\in I_{1}, m=0,\ldots,c-1\}$.
\end{prop}

\begin{figure}[ht]
\centering
\psfrag{0}{$0$}
\psfrag{1}{$1$}
\psfrag{x}{$x$}
\psfrag{I1}{$I_{1}$}
\psfrag{I2}{$I_{2}$}
\includegraphics{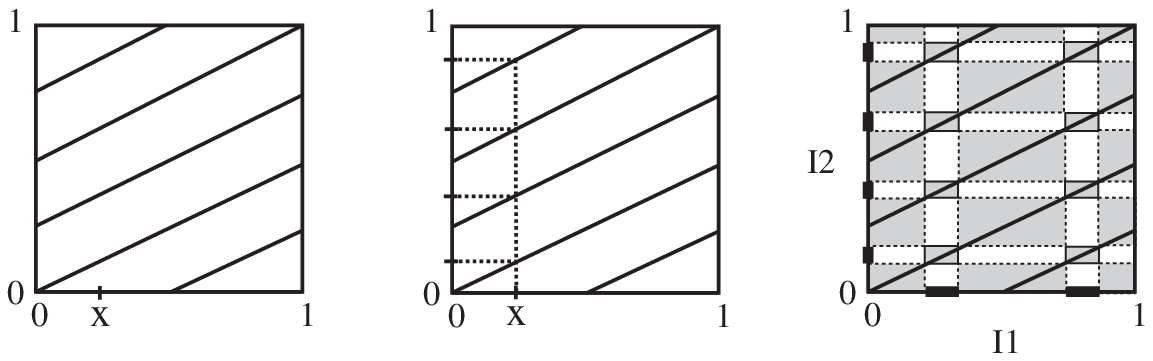}
\caption{We see $I_{1}$ with $|a|$-fold symmetry, $I_{2}$ with $|c|$-fold symmetry, and $\clorb\subset B$.  Also we see $x\in I_{1}$ and $\pi_2(\pi_1^{-1}(x)\cap\clorb)\subset I_2$.}
\label{fig:circlesontoruswithb}
\end{figure}

\begin{proof}
The ``if'' is clear; $\clorb=\clorb(0,0)=\{(x,y): ax+mb=cy,\,m=0,\ldots,g-1\}$ and by construction these lines lie entirely in $B=(I_{1}\times I_{2})\cup(I_{1}^{c}\times I_{2}^{c})$ (see Figure \ref{fig:circlesontoruswithb}).  So assume that the itineraries are the same.  Then $x\in I_1$ if and only if $\pi_2(\pi_1^{-1}(x)\cap\clorb)\subset I_2$, where $\pi_i:\T^2\to S^1$ is projection onto the $i$th coordinate. Since $\pi_2(\pi_1^{-1}(x)\cap\clorb)= \{ax/c+mb/c:x\in I_{1}, m=0,\ldots,|c|-1\}$, $I_2$ has $|c|$-fold symmetry (see Figure \ref{fig:circlesontoruswithb}).  The same argument, with 1 and 2 reversed, shows that $I_{1}$ has $|a|$-fold symmetry.
\end{proof}

\begin{proof}[Proof of Theorem \ref{thm:unique}]
Let $\alpha_{1},\alpha_{2}\in [0,1/2]\bs \Q$ and let $I_{1},I_{2}\subset S^{1}$ be finite unions of closed intervals  with no rotational symmetries.  Let $(a_{0},a_{1},a_{2},\ldots)$ and $(b_{0},b_{1},b_{2},\ldots)$ be the itineraries for $(\alpha_{1},I_{1})$ and $(\alpha_{2},I_{2})$, respectively.  

If $\alpha_{1}$ and $\alpha_{2}$ are not rationally related, then $\orb(0,0)$ is dense in $\T^{2}$ under $R_{\alpha_{1},\alpha_{2}}$.  In particular, $R^{n}(0,0)\not\in B$ for some $n\ge 0$.  In this case $a_{n}\ne b_{n}$.

Now suppose $\alpha_{1}$ and $\alpha_{2}$ are irrational but rationally related.  By Proposition \ref{prop:ratdep}, the itineraries are the same if and only if $a=1$, $c=1$ and $I_{1}=I_{2}$ (or  $a=1$, $c=-1$ and $I_{1}=-I_{2}$).  That is, $\alpha_{1}=b\pm\alpha_{2}$, or equivalently, $[\alpha_{1}]=[\alpha_{2}$].  Since $\alpha_{1},\alpha_{2}\in [0,1/2]$, $\alpha_1=\alpha_2$.
\end{proof}

It turns out that we can quantify the similarity of two itineraries.  That is, if we compare two itineraries term-by-term we can determine the asymptotic fraction of terms that are the same.  Let $\clorb=\clorb(0,0)$ and let $B=(I_{1}\times I_{2})\cup(I_{1}^{c}\times I_{2}^{c})$, as before.  Because $R_{\alpha_{1},\alpha_{2}}$ is uniquely ergodic when restricted to $\clorb$ (Proposition \ref{prop:torusorbits}), the fraction of the itineraries that are the same is precisely the fraction of $\clorb$ that intersects $B$ (\cite[Cor.~4.1.14]{KH}).  

If $\alpha_{1}$, $\alpha_{2}$, and $1$ are not rationally related then $\clorb=\T^{2}$ so $\clorb\cap B=B$.  If $\alpha_{1}$, $\alpha_{2}$, and 1 are rationally related, then $\clorb$ is  a finite set of circles, $k$, say.  Each circle corresponds to a line of rational slope $a/c$.  If $a/c$ is reduced, then the length of each circle is $\sqrt{a^{2} + c^{2}}$. (Recall that we use $l$ to denote the length of an interval (or a collection of intervals) in $S^1$, considered as $\R/\Z$.  Similarly, we will also use $l$ for Euclidean lengths in $\T^2$, considered as $\R^2/\Z^2$.)
Thus we have the following proposition. (Let $\delta_{a_i,b_i}=1$ if $a_i=b_i$ and 0 otherwise.)

\begin{prop}\label{prop:torusirrational}
Let $\alpha_{1},\alpha_{2}\in\R\bs\Q$, $(a_{0},a_{1},a_{2},\ldots)$ and $(b_{0},b_{1},b_{2},\ldots)$ be the itineraries for $(\alpha_{1},I_{1})$ and $(\alpha_{2},I_{2})$, and $\orb$ be the orbit of $0$.   
\begin{enumerate}
\item
If $1$, $\alpha_{1}$, and $\alpha_{2}$ are not rationally related, then \[\lim_{n\to\infty}\frac{1}{n}\sum_{i=0}^{n-1}\delta_{a_{i},b_{i}}=l(I_{1})l(I_{2})+l(I_{1}^{c})l(I_{2}^{c}).\]
\item
If $1$, $\alpha_{1}$, and $\alpha_{2}$ are rationally related and $\clorb$ is a set of $k$ circles of slope $a/c$ (reduced), then \[\lim_{n\to\infty}\frac{1}{n}\sum_{i=0}^{n-1}\delta_{a_{i},b_{i}}=\frac{l(B\cap \clorb)}{k\sqrt{a^{2}+c^{2}}}.\]
\end{enumerate}
\end{prop}

\section{Recovering $I$}
\label{sec:findI}

In the previous section we discussed the extent to which itineraries are unique.  We now address the question: given only an itinerary, is it possible to recover $[\alpha]$ and $I$?  First we point out that if we know $\alpha$ then it is relatively easy to recover $I$.  This follows from the fact that the orbit of $0$ is dense in $S^{1}$.  

\begin{thm}\label{thm:findI}
Let $r_{\alpha}:S^{1}\to S^{1}$ be a rotation by $\alpha\in\R\bs\Q$, $I\subset S^{1}$ be a finite union of closed intervals, and $(a_{0},a_{1},\ldots)$ be the itinerary of $0$.  Then $I=\cl\big(\{i\cdot\alpha:a_{i}=1\}\big)$.
\end{thm}
%\begin{thm}\label{beta}
%Let $r_{\alpha}$ be a rotation by $\alpha\in\R$ and let $I=[0,\beta)$ with $\beta\in(0,1)$.  Let $(a_{0},a_{1},\ldots)$ be the itinerary of a point $z\in S^{1}$.  If the itinerary is periodic with period $q$, then $\displaystyle\frac{p-1}{q}\le\beta< \frac{p}{q}$ where $\displaystyle p=\sum_{i=0}^{q-1}a_{i}$.  If the itinerary is not periodic, then $\displaystyle\beta=\lim_{n\to\infty}\frac{1}{n}\sum_{k=0}^{n-1}a_{k}$.
%\end{thm}

If we have only a finite portion $(a_{0},a_{1},\ldots,a_N)$ of the itinerary, then we can approximate $I$ as follows, assuming that there are both 1's and 0's in the partial itinerary.  For simplicity, assume that $I$ contains only one interval (the case where $I$ contains multiple intervals is similar).  
Rename the points $0,\alpha, 2\alpha,\ldots, N\alpha$ as $x_{0},x_{1},x_{2},\ldots,x_{N+1}=x_{0}$ so that they are  listed in counterclockwise order around the circle.  Let $x_L$ be the unique point such that $x_L$ is in $I$ and $x_{L-1}$ is not (this is determined by the corresponding values in the itinerary).  Similarly, let $x_R$ be the unique point such that $x_R$ is in $I$ and $x_{R+1}$ is not.  Then $I$ contains the interval $[x_L, x_R]$ and is contained in the interval $[x_{L-1}, x_{R+1}]$.  

The exact spacing between the points of the orbit depends in a complicated way on the Diophantine properties of $\alpha$ (see \cite{S2}).
However, we can estimate the spacing as follows.  Assume we have $\alpha =p/q \pm \ep$ (where $p/q$ is in lowest terms).  Since $q$ iterates of $r_{p/q}$ give $q$ points spaced $1/q$ apart, we see that the maximum spacing for the first $q$ iterates of $r_\alpha$, and thus the maximum error for the estimate of an endpoint of $I$, is $1/q + 2q\ep$.
 
Of course, if we did not know $\alpha$, but knew $\{\alpha\}$, then we could still find $I$.  It would be $I=\cl\big(\{i\cdot\{\alpha\}:a_{i}=1\}\big)$.  On the other hand, if we knew only $[\alpha]$, then there would be two possibilites for $I$, namely $I=\cl\big(\{i\cdot[\alpha]:a_{i}=1\}\big)$ and $I=\cl\big(\{-i\cdot[\alpha]:a_{i}=1\}\big)$. 

In fact, if we are only interested in finding the total length of $I$ we may use the fact that $r_{\alpha}$ is uniquely ergodic.  Again, recall from \cite[Cor.~4.1.14]{KH} that for a uniquely ergodic map the Birkhoff averages converge everywhere.  Thus we can find the size of $I$ by computing the asymptotic fraction of time the orbit of $0$ spends in $I$.  Of course, this is the same as determining the asymptotic fraction of terms in the itinerary that are 1.  That is, \[l(I)=\lim_{n\to\infty}\frac{1}{n}\sum_{k=0}^{n-1}a_{k}.\]
The exact rate of convergence depends, again, on the Diophantine properties of $\alpha$ (see \cite{Schoiss} and \cite[Thm.~4.6]{BI}), and may be arbitrarily slow.  Assuming that $\alpha =p/q \pm \ep$ (where $p/q$ is in lowest terms) and that $I$ is a single interval, we can use the fact that the spacing of the first $q$ iterates of $r_\alpha$ is $1/q \pm 2q\ep$ to conclude that $(S_q -1)(1/q - 2q\ep) \le l(I) \le (S_q +1)(1/q + 2q\ep)$, where $S_q= \sum_{k=0}^{q-1}a_{k}$.

\section{Recovering $[\alpha]$}
\label{sec:singleInt}

Finding $[\alpha]$ from only the itinerary is more difficult than finding $I$.  If we knew that after $n$ iterates the point had traveled around the circle $m=m(n)$ times, then $\displaystyle\frac{m+1}{n}\le[\alpha]<\frac{m}{n}.$ Since $\displaystyle\lim_{n\to\infty}\frac{m+1}{n}=\lim_{n\to\infty}\frac{m}{n}$, $\displaystyle[\alpha]=\lim_{n\to\infty}\frac{m(n)}{n}$. However, if we are given an itinerary for unknown $\alpha$ and $I$, it is not clear how to find $m(n)$ for a given $n$.

For the sake of simplicity we first assume $I$ is a single closed interval and we present the method for finding $[\alpha]$ in this case. In the next section we describe how to generalize this technique to the multi-interval case.

There is a case in which it is easy to find $[\alpha]$.  If the size of the rotation is less than or equal to the interval length ($[\alpha]\le l(I)$) and also less than the length of the complement of the interval ($[\alpha]<1-l(I))$, then the orbit will land in $I$ and in $S^{1}\bs I $ each time it goes around the circle.  In this case we could look at $(a_{0},a_{1},\ldots,a_{n-1})$ and count the number of blocks of 1's.  If there are $b=b(n)$ blocks then the orbit has gone around the circle either $m=b$, $m=b-1$, or $m=b-2$ times.
 Clearly $\displaystyle\frac mn \le [\alpha] \le \frac {m+1}n$, so we have $\displaystyle\frac {b-2}n \le [\alpha] \le \frac {b+1}n$.  Still, in order to use this technique we must be able to look at an itinerary and say with certainty that $[\alpha]<\min\{l(I),1-l(I)\}$.  Because the orbit of $0$ is dense and $I$ is closed, it is not difficult to see that $[\alpha]\le l(I)$ if somewhere in the itinerary there are two consecutive 1's.  Similarly, $[\alpha]<1-l(I)$ if and only if somewhere in the itinerary there are two consecutive 0's.  Thus we have the following lemma.

\begin{lemma}\label{lem:findalpha}
Let $r_{\alpha}:S^{1}\to S^{1}$ be a rotation by $\alpha\in\R\bs\Q$, $I\subset S^{1}$ be a single closed interval, and $(a_{0},a_{1},\ldots)$ be the itinerary of $0$.   If the itinerary  contains two consecutive 1's and two consecutive 0's then $\displaystyle \frac {b(n)-2}n \le [\alpha] \le \frac {b(n)+1}n$ where $b(n)$ denotes the number of blocks of 1's in $(a_{0},\ldots,a_{n-1})$. Thus $\displaystyle[\alpha]=\lim_{n\to\infty}\frac{b(n)}{n}$.
\end{lemma}

This technique fails when $l(I)<[\alpha]$, because the orbit may hop over the interval once, twice or many times before it lands in the interval again.  We encounter a similar problem when $1-l(I)<[\alpha]$. The key to handling these cases, as we shall see, is to work with a higher power of $r_{\alpha}$.  We have two lemmas, the first of which follows from the fact that $r^{k}_{\alpha}=r_{k\alpha}$.

\begin{lemma}\label{lem:kalpha}
If $(a_{0},a_{1},a_{2},\ldots)$ is the itinerary of $0$ under $r_{\alpha}$, then the itinerary of $0$ under $r_{k\alpha}$ is $(a_{0},a_{k},a_{2k},a_{3k},\ldots)$.
\end{lemma}

\begin{lemma}\label{lem:kalphasmall}
There is a $k>0$ such that the itinerary of $0$ under $r_{k\alpha}$ has two consecutive 0's and two consecutive 1's.
\end{lemma}
\begin{proof}
Since $\alpha$ is irrational, the values $0,[\alpha],[2\alpha],[3\alpha],\ldots$ are dense in $[0,\frac{1}{2}]$.  Thus, for some $k>0$, $[k\alpha]<\min\{l(I),1-l(I)\}$.  For this value of $k$ it is possible for the orbit of $0$ under $r_{k\alpha}$ to remain in $I$ for two consecutive iterates, and because the orbit is dense in $S^{1}$, it will occur.  Likewise, the orbit will eventually remain in $I^{c}$ for two consecutive iterates.
\end{proof}

We now describe the method for finding $[\alpha]$ when $[\alpha]\ge\min\{l(I),1-l(I)\}$.   By Lemma \ref{lem:kalphasmall} we can find a $k>0$ such that the itinerary for $0$ under $r_{k\alpha}$, $(a_{0},a_{k},a_{2k},a_{3k},\ldots)$, has two consecutive 0's and two consecutive 1's.  By Lemma \ref{lem:findalpha} this enables us to find $[k\alpha]$, which in turn gives us two possibilities for $\{k\alpha\}$: $[k\alpha]$ and $1-[k\alpha]$.  In either case we can use Theorem \ref{thm:findI} to find the corresponding intervals, $I_{1}$ and $I_{2}$.  Either $I_{1}=I$ or $I_{2}=I$.   Moreover, both of these possibilities for $\{k\alpha\}$ give $k$ possibilities for $[\alpha]$; namely, if $\{k\alpha\}=[k\alpha]$ then $\{\alpha\}=\{\frac{n}{k}+[k\alpha]\}$ for some $n=0,\ldots,k-1$ and if $\{k\alpha\}=1-[k\alpha]$ then $\{\alpha\}=\{\frac{n}{k}-[k\alpha]\}$ for some $n=0,\ldots,k-1$. 

Finally, we must test each of the $2k$ candidates to see which one is $\alpha$.   Given two candidates, $\alpha_1= \frac{n_1}{k}\pm[k\alpha]$ and $\alpha_2= \frac{n_2}{k}\pm[k\alpha]$, we solve algebraically for $N$ such that
$r_{\alpha_1}^N(0)\in I$ and $r_{\alpha_2}^N(0)\not\in I$.
We compare this to the original itinerary to eliminate one of the two candidates (if $a_N=0$, then $\alpha\ne\alpha_1$; if $a_N=1$, then $\alpha\ne\alpha_2$).  We continue eliminating candidates until we are down to one, which must be $\alpha$.

When $I$ has $p>1$ subintervals the situation is more difficult. Just as in the single interval case, finding $[\alpha]$ hinges on determining how many times an orbit has gone around the circle---if not for $r_{\alpha}$, then for a power of $r_{\alpha}$. In the single interval case, having two consecutive 1's and two consecutive 0's was sufficient to guarantee that the rotation was less than or equal to the length of both the interval and the complement of the interval.  In the multiple interval case we need to ensure that the rotation is less than the length of every subinterval of both $I$ and $I^{c}$. Thus the previous technique no longer works. 

If we know a value $k$ so that $[k\alpha]$ is smaller than every subinterval of both $I$ and $I^{c}$ \and we know $p$, the number of intervals in $I$,
 then we can proceed as before. Each block of 0's and each block of 1's in the itinerary $(a_{0},a_{k},a_{2k},\ldots)$ corresponds to a connected interval. Then $\displaystyle[k\alpha]=\lim_{n\to\infty}\frac{m(n)}{pn}$ where $m(n)$ denotes the number of blocks of 1's in $(a_{0},a_{k},\ldots,a_{k(n-1)})$. Then use the technique from Section \ref{sec:singleInt} to find $I$ and $\{\alpha\}$. 
The problem is that we do not know $k$ or $p$. So we proceed as follows.

Order the set $\{(k,p):k,p\in\Z^{+}\}$ to obtain a sequence $(k_{1},p_{1}),(k_{2},p_{2}),\ldots$ Use the technique described in the previous paragraph with $p=p_{1}$ and $k=k_{1}$. 
By uniqueness (Theorem \ref{thm:unique}), if we have the wrong $k$ and $p$, then there is no set of intervals that will give the correct itinerary under rotation by $[k\alpha]$.  Thus the $I$ that we construct will have more than $p$ components; at this point we recognize that our choice of $k$ and $p$ did not work and we try again with the next $(k_i,p_i)$.  However, we may have to examine an arbitrarily long piece of the itinerary to determine that $k$ and $p$ will not work. Thus, in the multiple interval case, we can never be certain in finite time that we have the right $\alpha$ and $I$, only that our estimates give the same itinerary for an arbitrary number of iterates.

Algorithms for determining whether a given sequence is the itinerary for a rotation by a given rational $p/q$ and a single interval are given in \cite{AGM} and \cite{STZ}.  No such recognition algorithm is possible in the current setting because
 {\em any} finite sequence can be part of an itinerary for a given irrational rotation for some collection of intervals.  We can, however, say something about the number of intervals necessary to produce a given itinerary.
 
 Let $p$ be the number of disjoint subintervals in $I$.  Let $B$ be the number of distinct lengths of maximal blocks in the itinerary consisting of all 1's, not counting the initial block if it starts with a 1.  (So for the sample itinerary on page~\pageref{sampleitinerary}, we have $B\ge2$.)  The three gap theorem, mentioned in section~\ref{sec:notaPrev}, implies that if $B>3$, then $p\ge2$.  We also have the following result.
 
\begin{prop}
Let $(a_{0},a_{1},a_{2},\ldots)$ be the itinerary of $0$ under $r_{\alpha}$.  Then either
	\begin{enumerate}
	\item every subinterval of $I^c$ has length $\le \alpha$, or
	\item $p \ge \log_4(B+1)$.
	\end{enumerate}
(The corresponding statement holds with the roles of $I$ and $I^c$ (and thus 0 and 1) reversed.)
\end{prop}

\begin{proof}
Assume that some subinterval $J$ of $I^c$ has length greater than $\alpha$.  Then an orbit of $r_\alpha$ cannot leave any subinterval of $I$ and return to the same subinterval without first passing through $J$, thus putting a 0 in the itinerary.  Since the three gap theorem says that the blocks of 1's corresponding to a single visit to a given subinterval $I_i$ of $I$ can take at most three values, say $a_i$, $b_i$, and $c_i$, the length of a block of 1's in the itinerary must have the form $\sum_{i=1}^p x_i$, where $x_i\in\{0,a_i,b_i,c_i\}$ and not all $x_i$'s are 0.  The number of possible lengths $B$ is thus bounded by $B\le (4^p-1)$, giving $p\ge \log_4(B+1)$.
\end{proof}

\section{Diffeomorphisms of the circle}
\label{sec:diffeos}

It is natural to return to the beginning of this article and ask all the same questions for orientation preserving homeomorphisms of the circle that we did for rigid rotations.   Recall that for such a homeomorphism $f:S^{1}\to S^{1}$ the \emph{rotation number} is \[\alpha=\alpha_{f}=\lim_{n\to\infty}\frac{F^{n}(x)-x}{n}\mod 1\] where $F:\R\to \R$ is any lift of $f$ and $x\in \R$.   (The limit always exists and is independent of  our choices of $F$ and $x$.)  In particular, we ask: if $f:S^{1}\to S^{1}$ is an orientation preserving homeomorphism with rotation number $\alpha\in\R$, $I\subset S^{1}$ is a finite union of closed intervals, and $(a_{0}, a_{1}, a_{2},\ldots)$ is the itinerary of $0$, can we recover $\alpha$ (or better yet, $f$ itself) and $I$?

Unfortunately, the answer is ``no.''  We encounter the analogous barriers to uniqueness that we did with rigid rotations (rational rotation numbers, clockwise versus counterclockwise ambiguity, and something like rotational symmetry).  Worse still, there exist homeomorphisms with irrational rotation numbers with orbits that are not dense.  In our case, that means that we could find a union of closed intervals $I$ that the orbit of 0 never intersects, so the itinerary of 0 is $(0,0,0,\ldots)$.  Even when we have dense orbits we cannot recover $I$ or $f$.

So, from here onward we require that our homeomorphisms are $C^{2}$ (actually, $C^{1}$ with bounded variation suffices), orientation preserving diffeomorphisms  with irrational rotation numbers.  The benefit of working in this setting is that we may take advantage of Denjoy's theorem (see \cite[Thm.\ 12.1.1]{KH}) which says that such a diffeomorphism with rotation number $\alpha$ is topologically conjugate to $r_{\alpha}$. (See \cite {STZ} for a discussion of the rational and non-smooth cases, and \cite{CGBSG} for a proof of the three-gap theorem for $C^2$ diffeomorphisms.)

In this case we can recover the rotation number of $f$, subject to the analogous restriction as for rotations.  We do so by turning the problem into one for a rigid rotation. 

Define $I_{\alpha}=\cl(\{i\cdot\alpha:a_{i}=1\}$.

\begin{prop}
\label{prop:diffeotorot}
Let $f$, $\alpha$, $I$, $(a_{0}, a_{1}, a_{2},\ldots)$, and $I_{\alpha}$ be defined as above.  Then $I_{\alpha}$ is a finite union of closed intervals and $(a_{0}, a_{1}, a_{2},\ldots)$ is the itinerary of $0$ for the rigid rotation $r_{\alpha}$ relative to $I_{\alpha}$.
\end{prop}
\begin{proof}  By Denjoy's theorem there exists a homeomorphism $\varphi:S^{1}\to S^{1}$ such that $f\circ\varphi=\varphi\circ r_{\alpha}$ and $\varphi(0)=0$.  Then $I_{\alpha}=\varphi(I)$, so it is clearly the finite union of closed intervals.  That $(a_{0}, a_{1}, a_{2},\ldots)$ is the itinerary of $0$ for $r_{\alpha}$ is clear from the definition of $I_{\alpha}$.
\end{proof}

The beauty of Proposition \ref{prop:diffeotorot} is that it allows us to turn questions about $f$ and $I$ into questions about $r_{\alpha}$ and $I_{\alpha}$.  Unfortunately, it also tells us that knowing the itinerary of 0 for $f$ cannot give us any information about $f$ except perhaps the rotation number, nor can it give us the exact location of the intervals $I$.  

Theorem \ref{thm:unique} gave us a uniqueness theorem for itineraries of rigid rotations.  It is possible to combine Theorem \ref{thm:unique} with Proposition \ref{prop:diffeotorot} to obtain an analogous uniqueness theorem for itineraries for diffeomorphisms (which we will not state).  Unfortunately, the symmetry condition is more difficult to verify, for it is the interval $I_{\alpha}$, not $I$, that must have no rotational symmetries.

If we know nothing about $\alpha$ or the symmetries of $I_{\alpha}$, then we proceed as follows.  Pretend the itinerary $(a_{0},a_{1},a_{2},\ldots)$ comes from a rigid rotation $r_{\alpha'}$ relative to the collection of intervals $I_{\alpha'}$ (with no symmetries).  The techniques of Sections \ref{sec:findI} and \ref{sec:singleInt} enable us to find $[\alpha']$ and $I_{\alpha'}$.  Then $\alpha=(m+[\alpha'])/n$ for some $m\in\Z$, $n\in\Z^{+}$ (implying that $I_{\alpha}$ has $n$-fold symmetry).  Thus we have the following theorem.

\begin{thm}
Let $f:S^{1}\to S^{1}$ be a $C^{2}$  orientation preserving diffeomorphism with rotation number $\alpha\in\R\bs\Q$ and $I\subset S^{1}$ be a finite union of closed intervals.  Let $\alpha'\in\R\bs\Q$ and $J\subset S^{1}$ be a finite union of closed intervals with no rotational symmetries.  If $(a_{0},a_{1},a_{2},\ldots)$ is the itinerary of 0 relative to $I$ under $f$ and of 0 relative to $J$ under $r_{\alpha'}$, then $\alpha=(m+[\alpha'])/n$ for some  $m\in\Z$, $n\in\Z^{+}$.
\end{thm}

We have the following interesting corollary to this theorem.

\begin{cor}
Let $f,g:S^{1}\to S^{1}$ be $C^{2}$  orientation preserving diffeomorphisms with rotation numbers $\alpha_{f},\alpha_{g}\in\R\bs\Q$.  Then $\alpha_{f}$, $\alpha_{g}$, and 1 are rationally related if and only if there exist finite unions of closed intervals $I_{f},I_{g}\subset S^{1}$ so that the itineraries of 0 under $f$ and $g$, relative to $I_f$ and $I_g$ respectively, are the same.
\end{cor}

\nocite{CN}
\bibliographystyle{plain} 
\bibliography{circleshift} 
\end{document}